\documentclass[12pt]{amsart}

\newcommand{\newsection}[1]{\setcounter{equation}{0} \section{#1}}
\newcommand{\bea}{\begin{eqnarray}}
\newcommand{\eea}{\end{eqnarray}}

\newcommand{\cle}{\mathcal{E}}
\newcommand{\clf}{\mathcal{F}}

\newcommand{\clh}{\mathcal{H}}

\newcommand{\cll}{\mathcal{L}}
\newcommand{\clm}{\mathcal{M}}

\newcommand{\clr}{\mathcal{R}}

\newcommand{\raro}{\rightarrow}

\newcommand{\cl}[1]{\mathcal{#1}}

\def \qed {\hfill \vrule height6pt width 6pt depth 0pt}
\def\textmatrix#1&#2\\#3&#4\\{\bigl({#1 \atop #3}\ {#2 \atop #4}\bigr)}
\def\dispmatrix#1&#2\\#3&#4\\{\left({#1 \atop #3}\ {#2 \atop #4}\right)}
\newcommand{\be}{\begin{equation}}
\newcommand{\ee}{\end{equation}}
\newcommand{\ben}{\begin{eqnarray*}}
\newcommand{\een}{\end{eqnarray*}}

\newcommand{\NI}{\noindent}

\newcommand{\bi}{\begin{itemize}}
\newcommand{\ei}{\end{itemize}}

\def\5{{5\superprime}}

\newtheorem{Theorem}{\sc Theorem}
\newtheorem{Lemma}{\sc Lemma}
\newtheorem{Proposition}{\sc Proposition}
\newtheorem{Corollary}{\sc Corollary}
\newtheorem{Definition}{\sc Definition}
\newtheorem{Example}{\sc Example}
\newtheorem{Remark}{\sc Remark}
\newtheorem{Note}{\sc Note}
\newtheorem{Question}{\sc Question}
\newtheorem{ass}{\sc Assumption}
\newcommand{\bt}{\begin{Theorem}}
\def\beginlem{\begin{Lemma}}
\def\beginprop{\begin{Proposition}}
\def\begincor{\begin{Corollary}}
\def\begindef{\begin{Definition}}
\def\beginexamp{\begin{Example}}
\def\beginrem{\begin{Remark}}
\def\beginq{\begin{Question}}
\def\beginass{\begin{ass}}
\def\beginnote{\begin{Note}}
\newcommand{\et}{\end{Theorem}}
\def\endlem{\end{Lemma}}
\def\endprop{\end{Proposition}}
\def\endcor{\end{Corollary}}
\def\enddef{\end{Definition}}
\def\endexamp{\end{Example}}
\def\endrem{\end{Remark}}
\def\endq{\end{Question}}
\def\endass{\end{ass}}
\def\endnote{\end{Note}}

\begin{document}

\title{On unitarily equivalent submodules}

\author[Douglas]{Ronald G. Douglas}
\author[Sarkar]{Jaydeb Sarkar}

\address{Texas A \& M University, College Station, Texas 77843, USA}

\email{rdouglas@math.tamu.edu, jsarkar@math.tamu.edu}



\begin{abstract}
The Hardy space on the unit ball in $\mathbb{C}^n$ provides examples of a quasi-free, finite rank Hilbert module which contains a pure submodule isometrically isomorphic to the module itself. For $n=1$ the submodule has finite codimension. In this note we show that this phenomenon can only occur for modules over domains in $\mathbb{C}$ and for finitely-connected domains only for Hardy-like spaces, the bundle shifts. Moreover, we show for essentially reductive modules that even when the codimension is infinite, the module is subnormal and again, on nice domains such as the unit ball, must be Hardy-like.
\end{abstract}

\maketitle

\newsection{Introduction}

\vspace{1.5cm}

One approach to multivariate operator theory is via the study of a Hilbert space 
which are modules over some natural algebra. Examples of such an algebra are $A(\Omega)$, which can be
defined for any bounded domains $\Omega$ in $\mathbb{C}^n$ as the 
completion, with respect to the supremum norm over $\Omega$, of the functions holomorphic on a 
neighborhood of the closure of $\Omega$. One way to obtain such Hilbert modules is as the 
closure of $A(\Omega)$ in $L^2(\mu)$ for some measure $\mu$ on the closure of 
$\Omega$. For volume measure, one obtains the Bergman space $L^2_a(\Omega)$ 
and, when $\partial\Omega$ is smooth, the Hardy space $H^2(\Omega)$ for 
surface measure on $\partial\Omega$. There is a class of modules called 
quasi-free to which these examples belong.

One consequence of the celebrated theorem of Beurling \cite{Beu} is that all non-zero 
submodules of $H^2(\mathbb{D})$, where $\mathbb{D}$ is the unit disk in 
$\mathbb{C}$, are unitarily equivalent to $H^2(\mathbb{D})$. For submodules of 
$H^2(\mathbb{D}^n)$ over $A(\mathbb{D}^n)$, $n>1$, some are unitarily equivalent to $H^2(\mathbb{D}^n)$ and 
some are not. For the Hardy module $H^2(\partial\mathbb{B}^n)$ over the ball algebra 
$A(\mathbb{B}^n)$, for the unit ball $\mathbb{B}^n$ in $\mathbb{C}^n$, the 
existence of inner functions on $\mathbb{B}^n$ \cite{Alek} established the existence of 
proper submodules of $H^2(\partial\mathbb{B}^n)$ that are unitarily equivalent to $H^2(\mathbb{B}^n)$. For the Bergman modules over the polydisk or the 
ball, one can show (cf. \cite{C-G,P,Rich}) that no proper submodule is unitarily equivalent to the 
Bergman module itself.

In this note we consider the question of which Hilbert modules have proper submodules 
unitarily equivalent to the original. If $U$ is an isometric module map on the Hilbert module
$\mathcal{M}$, then $U\mathcal{M}$ is a submodule of $\mathcal{M}$ unitarily 
equivalent to $\mathcal{M}$. Conversely, all such unitarily equivalent 
submodules have such a representation. Now $U$ is unitary iff $U\mathcal{M} = 
\mathcal{M}$. If $U$ is a nonconstant unitary module map on $\mathcal{M}$, although its existence yields some implications about the nature of $\mathcal{M}$, this is not the phenomenon we examine in this note. Here, we consider the case in which $\bigcap\limits^\infty_{k\ge 0} U^k 
\mathcal{M}=(0)$, which, for reasons which are apparent, we will call a 
\emph{pure} unitarily equivalent submodule. We are concerned with the question of when 
such submodules exist.  More precisely, what can one say about the algebra or the module in such a case. We prove three sets of results with some corollaries. 

First,  we show if $\dim \mathcal{M}\ominus U\mathcal{M}<\infty$, then $n=1$. Moreover,
for  a finite rank $k$, quasi-free Hilbert module $\mathcal{M}$ over 
$A(\mathbb{D})$, the existence of a pure unitarily equivalent submodule of finite codimension implies  
that $\mathcal{M}$ is unitarily equivalent to $H^2_{\mathcal{E}}(\mathbb{D})$ 
with $\dim \mathcal{E} = k$. If $\Omega$ is finitely connected with nice boundary, then the same result holds with bundle shifts (cf.\ \cite{A-D}) replacing the Hardy space. Second, we show that an essentially reductive, quasi-free 
Hilbert module $\mathcal{M}$ over $A(\Omega)$ for which there exists a pure unitarily 
equivalent submodule must be subnormal and for $\Omega= \mathbb{B}^n$, is unitarily equivalent to $H^2_{\mathcal{E}}(\mathbb{B}^n)$ with 
$\dim \mathcal{E}<\infty$.  All of the results lend support to the conjecture that isometrically isomorphic submodules of finite codimension  for which the original module is essentially reductive can occur only for \v Silov modules \cite{D-P}. Finally, for a class of measures $\mu$ on the closure of $\Omega$, we show that two submodules of $L^2_a(\mu)$ are isometrically isomorphic iff they are equal generalizing results of Richter \cite{Rich}, Putinar \cite{P}, and Guo--Hu--Xu \cite{GJX}.

An ancillary goal of this note is the development of techniques for the study of multivariate operator theory and we draw on methods from algebra and operator theory.

All the Hilbert spaces in this note are separable and are over
the complex field $\mathbb{C}$. For a Hilbert space $\clh$, we
denote the Banach space of all bounded linear operators by
$\cll(\clh)$.

We begin by recalling the definition of quasi-free Hilbert module over $A(\Omega)$ 
which was introduced in (\cite{DM1},\cite{DM2}) and is related to earlier ideas of Curto--Salinas \cite{C-S}.  The Hilbert space $\clm$ is 
said to be a contractive Hilbert module over $A(\Omega)$ if $\clm$ is a unital 
module over $A(\Omega)$ with module map $A(\Omega) \times \clm \raro \clm$ such 
that 
$$\|\varphi f \|_{\cl M} \leq \|\varphi\|_{A(\Omega)} \|f\|_{\cl M}$$
 for $\varphi$ in 
$A(\Omega)$ and $f$ in $\clm$.

A Hilbert space $\clr$ is said to be a quasi-free Hilbert module of rank $m$ 
over $A(\Omega)$, $1 \leq m \leq \infty$, if it is obtained as the completion of 
the algebraic tensor product $A(\Omega) \otimes \ell^2_m$ relative to an inner 
product such that:

(1) $\mbox{eval}_{\pmb{z}}\colon\ A(\Omega) \otimes l^2_m \raro l^2_m$ is bounded for $\pmb{z}$ in 
$\Omega$ and locally uniformly bounded on $\Omega$;

(2) $\|\varphi(\sum \theta_i \otimes x_i)\|_{\clr} = \|\sum \varphi \theta_i \otimes 
x_i\|_{\clr} \leq \|\varphi\|_{A(\Omega)} \|\theta_i \otimes x_i\|_{\clr}$ for 
$\varphi$, $\{\theta_i\}$ in $A(\Omega)$ and $\{x_i\}$ in $\ell^2_m$; and 

(3) For $\{F_i\}$ a sequence in $A(\Omega) \otimes \ell^2_m$ which is Cauchy in 
the $\clr$-norm, it follows that $\mbox{eval}_{\pmb{z}}(F_i) \raro 0$ for all $\pmb{z}$ in 
$\Omega$ if and only if $\|F_i\|_{\clr} \raro 0$. 

If $I_{\pmb{\omega}_0}$ denotes the maximal ideal of polynomials in $\mathbb{C}[\pmb{z}] = \mathbb{C}[z_1,\ldots, z_n]$ which  vanish at $\pmb{\omega}_0$ for some $\pmb{\omega}_0$ in $\Omega$, then the Hilbert module $\mathcal{M}$ is said to be semi-Fredholm at $\pmb{\omega}_0$ if $\dim \mathcal{M}/{I_{\pmb{\omega}_0}\cdot \mathcal{M}} = m$ is finite (cf.\ \cite{E}). In particular, note that $\mathcal{M}$ semi-Fredholm at $\pmb{\omega}_0$ implies that $I_{\pmb{\omega}_0}\mathcal{M}$ is a closed submodule of $\mathcal{M}$.

One can show  that $\pmb{\omega}\to \mathcal{R}/I_{\pmb{\omega}}\cdot\mathcal{R}$ can be made into a rank $m$ Hermitian holomorphic vector bundle over $\Omega$ if $\mathcal{R}$ is semi-Fredholm at $\pmb{\omega}$ in $\Omega$, $\dim \mathcal{R}/{I_{\pmb{\omega}}\cdot \mathcal{R}}$ is constant $m$ and $\mathcal{R}$ is quasi-free, $1\le m<\infty$. Actually, all we need here is that the bundle obtained is real-analytic which is established in
(\cite{C-S}, Theorem 2.2).

A quasi-free Hilbert module of rank $m$ is a reproducing kernel Hilbert space 
with the kernel 
$$K(\pmb{w}, \pmb{z}) = \mbox{eval}_{\pmb{w}} \mbox{eval}_{\pmb{z}}^* \colon\ \Omega \times \Omega 
\raro \cll(\ell^2_m).$$ 

Before continuing to the existence question of isometrically isomorphic submodules, let us consider a useful property of quasi-free Hilbert modules. For $\varphi$ in $A(\Omega)$, let $M_\varphi$ denote the bounded operator on the Hilbert module defined by module multiplication by $\varphi$.

\vspace{0.3in}

\begin{Lemma}\label{lem1}
If $\mathcal{R}$ is a quasi-free Hilbert module over $A(\Omega)$ which is semi-Fredholm for $\pmb{\omega}$ in $\Omega$, then $\bigcap\limits_{\varphi\in I_{\pmb{\omega}}} \ker M_\varphi = (0)$ for $\pmb{\omega}$ in $\Omega$.
\end{Lemma}

\NI\textbf{Proof.} First, recall that the projection-valued function $P_{\pmb{\omega}}$ is real-analytic (cf.\ \cite{C-S}), where $P_{\pmb{\omega}}$ is the projection onto the closed submodule $I_{\pmb{\omega}}\cdot \mathcal{R}$ of $\mathcal{R}$.

Second, for $f$ in $\bigcap\limits_{\varphi\in I_{\pmb{\omega}}} \ker M_\varphi$, we have $M_\varphi f = \varphi(\pmb{\omega})f$ for all $\varphi$ in $A(\Omega)$. Let $\tilde{\pmb{\omega}}$ be in $\Omega\backslash\{\pmb{\omega}\}$ and $\omega_i \ne \tilde\omega_i$ for some $1\le i \le m$. Then $(M_{z_i}-\tilde\omega_i)f = (\omega_i-\tilde\omega_i)f$ and thus $f$ is in $I_{\tilde{\pmb{\omega}}}\cdot \mathcal{R}$ for $\tilde{\pmb{\omega}}\ne \pmb{\omega}$. Thus $(I-P_{\tilde{\pmb{\omega}}})f = 0$ for $\pmb{\omega} \ne \tilde{\pmb{\omega}}$ which implies since $I-P_{\pmb{\omega}}$ is real-analytic that $(I-P_{\pmb{\omega}})f=0$ or $f$ is in $I_{\pmb{\omega}}\cdot \mathcal{R}$ also. However, $\bigcap\limits_{\varphi\in I_{\pmb{\omega}}} \ker M_\varphi \subseteq \bigcap\limits_{\tilde{\pmb{\omega}}\in\Omega} I_{\tilde{\pmb{\omega}}}\cdot \mathcal{R} = (0)$ which completes the proof.\qed

\vspace{0.2cm}

Note the above proof does not require that $\mathbb{C}[\pmb{z}]$ is dense in $A(\Omega)$.

\section{Finite Codimension Case}

\indent 

Many of the questions concerning pure unitarily equivalent submodules can be reduced to questions about Toeplitz operators as follows.

Recall that if $P$ is the Szeg\"o projection from 
$L^2_{\mathcal{E}}(\mathbb{T})$ to $H^2_{\mathcal{E}}(\mathbb{D})$, and 
$\varphi$ is a function in $L^\infty_{\mathcal{L}(\mathcal{E})}(\mathbb{T})$, then 
the Toeplitz operator $T_\varphi$ is defined on $H^2_{\mathcal{E}}(\mathbb{D})$ 
so that $T_\varphi f = P(\varphi \cdot f)$, where $\varphi\cdot f$ denotes the 
function defined by $(\varphi\cdot f)(e^{it}) = \varphi(e^{it})f(e^{it})$ for 
$e^{it}$ in $\mathbb{T}$.

We begin with a well-known lemma.

\vspace{0.3in}

\begin{Lemma}\label{Toep}
If $\varphi$ is a function in $H^{\infty}_{\cll(\cle)}(\mathbb{D})$ such that 
the Toeplitz operator $T_{\varphi}$ is bounded below on $H^2_{\cle}$, then the 
Laurent multiplication operator $L_{\varphi}$ is bounded below on 
$L^2_{\cle}(\mathbb{T})$.
\end{Lemma}

\NI\textbf{Proof.} Observe that the set $\{e^{-iNt} f \colon\ N \in \mathbb{Z}_{+}, f 
\in H^2_{\cle}(\mathbb{D})\}$ is dense in $L^2_{\cle}(\mathbb{T})$ and for each 
$e^{-iNt} f$ in that dense set, we have
$$\|L_{\varphi}(e^{-iNt} f)\| = \|e^{-iNt} T_{\varphi} f \| = |e^{-iNt}| 
\|T_{\varphi} f\| \geq \epsilon \|e^{-iNt} f \|,$$ 
where $\epsilon >0$ is the lower bound for $T_\varphi$. \qed

\vspace{0.3in}

Let us now consider how Toeplitz operators (cf. \cite{DB}) enter the picture. Let $\mathcal{M}$ be a Hilbert module over $A(\Omega)$ with $U$ an isometric module map which satisfies $\bigcap\limits^\infty_{k=0} U^k \mathcal{M}  =(0)$, that is, $U\mathcal{M}$ is pure. If we set $\mathcal{E} = \mathcal{M} \ominus U\mathcal{M}$, then there exists a canonical isomorphism $\Psi\colon \ H^2_{\mathcal{E}} (\mathbb{D})\to \mathcal{M}$ such that $\Psi T_z = U\Psi$. If $M_{z_i}$ denotes the operator on $\mathcal{M}$ defined by module multiplication by the coordinate function $z_i$, then  $X_i = \Psi^*M_{z_i}\Psi$ is an operator on $H^2_{\mathcal{E}}(\mathbb{D})$ which commutes with $T_z$. Hence, there exists a function $\varphi_i$ in $H^\infty_{\mathcal{L}(\mathcal{E})}(\mathbb{D})$ such that $X_i = T_{\varphi_i}$. Moreover, since the $\{M_{z_i}\}$ commute, so do the $\{X_i\}$ and hence the functions $\{\varphi_i\}$ commute pointwise a.e.\ on $\mathbb{T}$. In particular, $\mathcal{R}$ is isomorphic as a module  over $\mathbb{C}[\pmb{z}]$ to that obtained by letting $z_i$ act on $H^2_{\mathcal{E}}(\mathbb{D})$ by $T_{\varphi_i}$. Hence, the basic question is what Hilbert modules can be so represented in this form.

We summarize this construction as follows:
\vspace{0.4in}

\begin{Proposition}
Let $\mathcal{M}$ be a Hilbert module over $A(\Omega)$ for which there exists an isometric module map $U$ satisfying $\bigcap\limits^\infty_{k=0} U^k\mathcal{M} = (0)$. Then there exists an isomorphism $\Psi\colon \ H^2_{\mathcal{E}}(\mathbb{D})\to \mathcal{M}$ with $\mathcal{E} = \mathcal{M}\ominus U\mathcal{M}$ and a commuting $n$-tuple of functions $\{\varphi_i\}$ in $H^\infty_{\mathcal{L}(\mathcal{E})}(\mathbb{D})$ so that $U = \Psi T_z\Psi^*$ and $M_{z_i} = \Psi T_{\varphi_i}\Psi^*$ for $i=1,2,\ldots, n$.
\end{Proposition}

If we are to reach conclusions about $\Omega$, then we must find a closer connection between the Hilbert module $\mathcal{R}$ and $\Omega$. One possibility is to assume something about the Hilbert--Samuel polynomial $h^{\pmb{\omega}_0}_{\mathcal{R}}$ in $\mathbb{C}[z]$ for $\mathcal{R}$ (cf. \cite{A1}, \cite{D-Y}, \cite{DE}). Recall that $h^{\pmb{\omega}_0}_{\mathcal{R}}$ is a polynomial in one variable for which $h^{\pmb{\omega}_0}_{\mathcal{R}}(k) = \dim_{\mathbb{C}} I^k_{\pmb{\omega}_0}\cdot \mathcal{R}/I^{k+1}_{\pmb{\omega}_0}\cdot \mathcal{R}$ for all $k\ge N_{\mathcal{R}}$ for some positive integer $N_{\mathcal{R}}$. Here we are assuming that $\mathcal{R}$ is semi-Fredholm at $\pmb{\omega}_0$.

Consider rank $k$ quasi-free Hilbert modules $\mathcal{R}$ and $\widetilde{\mathcal{R}}$ over $A(\Omega)$ with $1\le k<\infty$. Following Lemma $1$ of \cite{DM2}, construct the rank $k$ quasi-free Hilbert module $\Delta$, which is the graph of a closed densely defined module map from $\mathcal{R}$ to $\widetilde{\mathcal{R}}$ obtained as the closure of the set $\{\varphi f_i\oplus \varphi g_i\colon\ \varphi\in A(\Omega)\}$, where $\{f_i\}$ and $\{g_i\}$ are generators for $\mathcal{R}$ and $\widetilde{\mathcal{R}}$, respectively. Then the module map $X\colon \ \Delta\to \mathcal{R}$ defined by $f_i\oplus g_i \to f_i$ is bounded, one-to-one and has dense range.

If we consider the adjoint $X^*\colon \ \mathcal{R}\to \Delta$, then for fixed $\pmb{\omega}_0$ in $\Omega$, $X^*(I_{\pmb{\omega}_0}\cdot \mathcal{R})^\bot \subset (I_{\pmb{\omega}_0}\cdot \Delta)^\bot$. Since the rank of $\Delta$ is also $k$, this map is an isomorphism. Let $\{\gamma_i(\pmb{\omega})\}$ be anti-holomorphic functions from a neighborhood $\Omega_0$ of $\pmb{\omega}_0$ to $\mathcal{R}$ such that $\{\gamma_i(\pmb{\omega})\}$ spans $(I_{\pmb{\omega}} \cdot \mathcal{R})^\bot$ for $\pmb{\omega}$ in $\Omega_0$. Then $\left\{\frac{\partial^{\pmb{\alpha}}}{\partial\pmb{z}^{\pmb{\alpha}}} \gamma_i(\pmb{\omega})\right\}_{|\pmb{\alpha}|< k}$ forms a basis for $(I^k_{\pmb{\omega}}\cdot \mathcal{R})^\bot$ for $k=0,1,2,\ldots$ using the same argument as in Section 4 in \cite{C-S} and Section 4 in \cite{DMV}. Similarly, since $\{X^*\gamma_i(\pmb{\omega})\}$ is a basis for $(I_{\pmb{\omega}}\cdot \Delta)^\bot$, we see that $X^*$ takes $(I^k_{\pmb{\omega}}\cdot \mathcal{R})^\bot$ onto $(I^k_{\pmb{\omega}}\cdot \Delta)^\bot$ for $k=0,1,2,\ldots$~. Therefore, $\dim(I^k_{\pmb{\omega}}\cdot \mathcal{R})^\bot = \dim(I^k_{\pmb{\omega}}\cdot\Delta)^\bot$ for all $k$.
Hence, $h^{\mathcal{R}}_{\pmb{\omega}} = h^\Delta_{\pmb{\omega}}$ for all $\pmb{\omega}$ in $\Omega$. Interchanging the roles of $\mathcal{R}$ and $\widetilde{\mathcal{R}}$ we have established the following result.

\vspace{0.4in}
\begin{Lemma}
If $\mathcal{R}$ and $\widetilde{\mathcal{R}}$ are semi-Fredholm, having the same finite rank, quasi-free Hilbert modules over $A(\Omega)$, then $h^{\mathcal{R}}_{\pmb{\omega}} \equiv h^{\widetilde{\mathcal{R}}}_{\pmb{\omega}}$ for $\pmb{\omega}$ in $\Omega$.
\end{Lemma}

In particular, one can calculate the Hilbert--Samuel polynomial by considering only the Bergman module over $A(\Omega)$ since $h^{\mathcal{R}\otimes \mathbb{C}^k}_{\pmb{\omega}}\equiv kh^{\mathcal{R}}_{\pmb{\omega}}$ for all finite integers $k$. To accomplish that we can reduce to the case of a ball as follows.

Let $\Omega$ be a bounded domain in $\mathbb{C}^n$ and $B_\varepsilon(\pmb{\omega}_0)$ be a ball with radius $\varepsilon$ centered at $\pmb{\omega}_0$, whose  closure is contained in $\Omega$. An easy argument shows that the map $X\colon \ L^2_a(\Omega)\to L^2_a(B_\varepsilon(\pmb{\omega}_0))$ defined by $Xf\equiv f|_{B_\varepsilon(\pmb{\omega}_0)}$ for $f$ in $L^2_a(\Omega)$ is bounded (actually compact), one-to-one and has dense range. Moreover, one can repeat the above argument for $\pmb{\omega}$ in $B_\varepsilon(\pmb{\omega}_0)$ to conclude that $h^{L^2_a(\Omega)}_{\pmb{\omega_0}} \equiv h^{B_\varepsilon(\pmb{\omega}_0)}_{\pmb{\omega_0}}$. The proof is completed by considering the Hilbert--Samuel polynomials at $\pmb{\omega}_0$ of the Bergman module for the ball $B_\varepsilon(\pmb{\omega}_0)$ for some $\varepsilon>0$ which is  centered   at $\pmb{\omega}_0$. This calculation reduces to that of the module $\mathbb{C}[\pmb{z}]$ over the algebra $\mathbb{C}[\pmb{z}]$ since the monomials in $L^2_a(B_\varepsilon(\pmb{\omega}_0))$ are orthogonal. Hence, $h^{L^2_a(B_\varepsilon(\pmb{\omega}_0)}_{\pmb{\omega}_0}(k) = \binom{n+k-1}n$ and we obtain:

\vspace{0.4in}

\begin{Proposition}\label{pro2}
If $\mathcal{R}$ is a quasi-free Hilbert module over $A(\Omega)$ for $\Omega\subset \mathbb{C}^n$ which is semi-Fredholm for $\pmb{\omega}$ in a neighborhood of $\pmb{\omega}_0$ in $\Omega$ with constant codimension, then $h^{\omega_0}_{\mathcal{R}}$ has degree  $n$.
\end{Proposition}

On the other hand, if there exists a pure isometrically isomorphic submodule of finite codimension, the Hilbert--Samuel polynomial is linear.

\begin{Proposition}\label{pro3}
If $\mathcal{M}$ is semi-Fredholm at $\pmb{\omega}_0$ in $\Omega$ and $\mathcal{N}$ is a pure isometrically isomorphic submodule of $\mathcal{M}$ having finite codimension in $\mathcal{M}$, then $h^{\pmb{\omega}_0}_{\mathcal{M}}$ has degree at most one.
\end{Proposition}

\NI\textbf{Proof.} As before, the existence of $\mathcal{N}$ in $\mathcal{M}$ yields a module isomorphism of $\mathcal{M}$ with $H^2_{\mathcal{E}}(\mathbb{D})$ for $\mathcal{E} = \mathcal{M}\ominus \mathcal{N}$. Assume that $\pmb{\omega}_0 = \pmb{0}$ for simplicity and note that the assumption that $\mathcal{M}$ is semi-Fredholm at $\pmb{\omega}_0 = \pmb{0}$ implies that $M_{z_1}\cdot \mathcal{M} +\cdots+ M_{z_m}\cdot \mathcal{M}$ has finite codimension in $\mathcal{M}$. Hence, $\widetilde{\mathcal{N}} = T_{\varphi_1}\cdot H^2_{\mathcal{E}}(\mathbb{D}) +\cdots+ T_{\varphi_m}\cdot H^2_{\mathcal{E}}(\mathbb{D})$ has finite codimension in $H^2_{\mathcal{E}}$. Moreover, $\widetilde{\mathcal{N}}$ is invariant under the action of $T_z$. Therefore, by the Beurling--Lax--Halmos Theorem (cf.\ \cite{S-N}), there is an inner function $\Theta$ in $H^\infty_{\mathcal{L}(\mathcal{E})}(\mathbb{D})$ for which $\widetilde{\mathcal{N}} = \Theta H^2_{\mathcal{E}}(\mathbb{D})$. Further, since $\widetilde{\mathcal{N}}$ has finite codimension in $H^2_{\mathcal{E}}(\mathbb{D})$ and the dimension of $\mathcal{E}$ is finite, it follows that the matrix  entries of $\Theta$ are rational functions with poles outside the closed unit disk and $\Theta(e^{it})$ is unitary for $e^{it}$ in $\mathbb{T}$ (cf.\ \cite{S-N}, Chapter VI, Section 4).

Now the determinant, $\det \Theta$, is a scalar-valued rational inner function in $H^\infty(\mathbb{D})$ and hence is a finite Blaschke product. Using Cramer's Rule one can show that $(\det\Theta)H^2_{\mathcal{E}}(\mathbb{D}) \subseteq\Theta H^2_{\mathcal{E}}(\mathbb{D})$ (cf. \cite{He}, Theorem 11) which implies that
\[
\dim_{\mathbb{C}} H^2_{\mathcal{E}}(\mathbb{D})/\Theta H^2_{\mathcal{E}}(\mathbb{D}) \le \dim_{\mathbb{C}} H^2_{\mathcal{E}}(\mathbb{D})/(\det \Theta)H^2_{\mathcal{E}}(\mathbb{D}).
\]
Continuing, we have 
\begin{align*}
\Psi(I^2_{\pmb{\omega}_0}\cdot \mathcal{M}) &= \Psi\left(\bigvee^n_{i,j=1} X_iX_j \mathcal{M}\right)\\
&= \bigvee^n_{i,j=1} T_{\varphi_i}T_{\varphi_j} H^2_{\mathcal{E}}(\mathbb{D})\\
&= \bigvee^n_{i=1} T_{\varphi_i}(\Theta H^2_{\mathcal{E}}(\mathbb{D})) \supseteq \bigvee^n_{i=1} T_{\varphi_i}(\det \Theta)H^2_{\mathcal{E}}(\mathbb{D})\\
&\supseteq \bigvee^n_{i=1} \det\Theta(T_{\varphi_i} H^2_{\mathcal{E}}(\mathbb{D})) = (\det\Theta) \Theta H^2_{\mathcal{E}}(\mathbb{D})\\
&\supseteq (\det\Theta)^2 H^2_{\mathcal{E}}(\mathbb{D}).
\end{align*}
Therefore, we have
\[
\dim(\mathcal{M}/I^2_{\pmb{\omega}_0}\cdot \mathcal{M}) \le \dim H^2_{\mathcal{E}}(\mathbb{D})/(\det \Theta)^2 H^2_{\mathcal{E}}(\mathbb{D}).
\]
Proceeding by induction, we obtain for each positive integer $k$ that
\[
\dim(\mathcal{M}/I^k_{\pmb{\omega}_0}\cdot \mathcal{M}) \le \dim H^2_{\mathcal{E}}(\mathbb{D})/(\det \Theta)^k H^2_{\mathcal{E}}(\mathbb{D}).
\]
If $D$ is the dimension  of $H^2/(\det \Theta)H^2$, then we have
\[
h^{\pmb{\omega}_0}_{\mathcal{M}}(k) \le \dim H^2_{\mathcal{E}}(\mathbb{D})/(\det\Theta)^k H^2_{\mathcal{E}}(\mathbb{D}) = k D\dim \mathcal{E} \quad \text{for}\quad k\ge N_{\mathcal{R}}.
\]
Hence, the degree of $h^{\pmb{\omega}_0}_{\mathcal{M}}$ is at most one.\qed

\vspace{0.3in}

Observe that $\mathcal{M}$ is not required to be quasi-free in this proposition. This proof extends one given by Fang in (\cite{XF}, Proposition 29). In particular, he shows that a necessary condition that one can represent a commuting $n$-tuple of isometries using inner functions on $H^2_{\mathcal{E}}(\mathbb{D})$, with $\dim \mathcal{E}<\infty$, is that the Hilbert--Samuel polynomial for the corresponding Hilbert module over $\mathbb{C}[\pmb{z}]$ is linear. One might predict a generalization of this result to show that the analogous criteria holds for inner functions in $H^{\infty}_{\cll(\mathcal{E})}(\mathbb{D}^k)$, $\dim \mathcal{E}<\infty$, with the degree of the corresponding Hilbert--Samuel polynomial being less than or equal to $k$. However, there is a gap  in using the preceding argument since the Beurling-Lax-Halmos Theorem is false for $H^2_{\cle}(\mathbb{D}^k)$ for $k>1$.

We combine Proposition 2 and 3 to obtain our first main result.

\vspace{0.4in}
\begin{Theorem}\label{thm2}
If $\mathcal{R}$ is a semi-Fredholm, quasi-free Hilbert module over $A(\Omega)$ with $\Omega\subset \mathbb{C}^n$ having a pure isometrically isomorphic submodule of finite codimension, then $n=1$.
\end{Theorem}

\vspace{0.4cm}

Hence one can immediately reduce to the case of domains $\Omega$ in $\mathbb{C}$ if there exists a pure isometrically isomorphic submodule of finite codimension.

Now let us consider what conclusions we can draw if we assume that $\Omega = \mathbb{D}$. In this case, $M_z$, module multiplication by $z$ on $\mathcal{R}$, corresponds to $T_\varphi$ on $H^2_{\mathcal{E}}(\mathbb{D})$ for some $\varphi$ in $H^\infty_{\mathcal{L}(\mathcal{E})}(\mathbb{D})$.

\vspace{0.4in}

\begin{Theorem}\label{thm3}
Let $\clm$ be a finite rank, quasi-free Hilbert module over $A(\mathbb{D})$ 
which is semi-Fredholm for $\omega$ in $\mathbb{D}$. Assume there exists a pure module isometry $U$ such that 
$\dim \clm/ U \clm <~\infty$. Then $\clm$ and  $H^2(\mathbb{D})$ are  
$A(\mathbb{D})$-module isomorphic.
\end{Theorem}

\NI\textbf{Proof.} 
As in Proposition 1, we can assume that $\clm 
\cong H^2_{\cle}(\mathbb{D})$, where $\cle = \clm \ominus U \clm$ with $\dim \cle 
< \infty$ and $U$ corresponds to $T_z$. Let $T_{\varphi}$ denote the operator on $H^2_{\mathcal{E}}(\mathbb{D})$ unitarily equivalent to  module multiplication 
by $z$ on $\mathcal{M}$, where $\varphi$ is in 
$H^{\infty}_{\cll(\cle)}(\mathbb{D})$ with $\| \varphi (z) \| \leq 1$ for all $z$ 
in $\mathbb{D}$.

Since the operator $T_{\varphi}$ is defined by module multiplication on 
$H^2_{\cle}(\mathbb{D})$ and the corresponding $A(\mathbb{D})$-module has finite rank,  it is enough to show that $\varphi$ is inner. Hence $T_\varphi$ would be a 
pure isometry so that $H^2_{\cle}(\mathbb{D})$ and $H^2(\mathbb{D})$ would be
$A(\mathbb{D})$-module isomorphic.

Since the range of $T_\varphi-\omega I$ has finite codimension in $H^2_{\mathcal{E}}(\mathbb{D})$ for $\omega$ in $\mathbb{D}$, it follows that the operator $T_{\varphi} - w I$ has closed range 
for each $w$ in $\mathbb{D}$. Now  $\mbox{Ker} (T_{\varphi} - w I) = 
\{0\}$ by Lemma \ref{lem1} so that $T_{\varphi} - w I$ is bounded below. Then by Lemma \ref{Toep},  
$(L_{\varphi} - w I)$ is bounded below on $L^2_{\cle}(\mathbb{D})$.

For each $w$ in $\mathbb{D}$ and $k$ in $\mathbb{N}$ let us consider the set 
$$E^{w}_k = \{ e^{it}\colon\  \mbox{dist}(\sigma(\varphi(e^{it})), w) < \frac{1}{k}\},$$
where $\sigma(\varphi(e^{it}))$ denotes the spectrum of the matrix  $\varphi(e^{it})$.

Then either $\sigma(\varphi(e^{it})) \subset \mathbb{T}$ a.e or there exists a 
$w_0$ in $\mathbb{D}$ such that  $m(E^{w_0}_k) > 0$ for all $k \in \mathbb{N}$. 
In the latter case, we can find a sequence of functions $\{f_k\}$ in $L^2_{\mathcal{E}}(\mathbb{T})$ such that 
$f_k$ is supported on $E^{w_0}_k$,
$\|f_k(e^{it})\| = 1$ for $e^{it}$ in $E^{\omega_0}_k$
 and 
$$\| \varphi(e^{it}) f_k (e^{it}) - w_0 
f_k (e^{it}) \| \leq \frac{1}{k}.$$
Thus we obtain that 
$$\| (L_{\varphi} - w_0) f_k \| \leq \frac{1}{k}\|f_k\|$$
 for all $k$ in $\mathbb{N}$, which contradicts the fact that $L_{\varphi} - 
w_0I$ is bounded below. Hence, $\sigma(\varphi(e^{it})) \subset \mathbb{T}$,
a.e.\ and hence $\varphi(e^{it})$ is unitary a.e. Therefore, $T_{\varphi}$ is a pure isometry and the Hilbert module $H^2_{\cle}(\mathbb{D})$ determined by $T_\varphi$ is $A(\mathbb{D})$-module isomorphic with $H^2(\mathbb{D})$.\qed

\vspace{0.3in}

If one attempts to extend the above proof to the case in which $\dim \mathcal{M}/U\mathcal{M}$ is infinite, one must confront the fact that there exist non-unitary contraction operators with spectrum containe in $\partial\mathbb{D}$. However, it  seems possible that the result is still valid in this case.

This result can't be extended to the case in which $U$ is not pure. For example, for $\mathcal{M} = H^2(\mathbb{D}) \oplus L^2_a(\mathbb{D})$, one could take $U = M_z\oplus I$.

The above proof can be extended to the case of a finitely-connected domain $\Omega$ with a nice boundary, that is, $\Omega$ for which $\partial\Omega$ is the finite union of simple closed curves. First, we must recall the notion of the bundle shift $H^2_\alpha(\Omega)$ for $\Omega$ determined by the unitary representation $\alpha$ of the fundamental group $\pi_1(\Omega)$ of $\Omega$. Then $H^2_\alpha(\Omega)$ is the Hardy space of holomorphic sections of the flat unitary bundle over $\Omega$ determined by $\alpha$ (cf.\ \cite{A-D}).

\vspace{0.4in}
\begin{Theorem}\label{thm4}
Let $\mathcal{M}$ be a finite rank, quasi-free Hilbert module over $A(\Omega)$, where $\Omega$ is a finitely-connected domain in $\mathbb{C}$ with  nice boundary, which is semi-Fredholm for $\omega$ in $\Omega$. Let $U$ be a pure module isometry such that $\dim \mathcal{M}/U\mathcal{M} < \infty$. Then there is a unitary representation of $\pi^1(\Omega)$ on some finite dimensional Hilbert space such that $\mathcal{M}$ and the bundle shift $H^2_\alpha(\Omega)$, are $A(\Omega)$-module isomorphic.
\end{Theorem}

\NI\textbf{Proof.} If we proceed as in the proof of the previous theorem, then the operators $\varphi(e^{it})$ on $\mathcal{E}$ for $e^{it}$ in $\partial\Omega$ have  $\text{clos } \Omega$ as a spectral set and the analogous argument shows that the eigenvalues of $\varphi(e^{it})$ lie in $\partial\Omega$. As a consequence $\varphi(e^{it})$ is normal. Hence multiplication $L_\varphi$ by $\varphi$ on $L^2_{\mathcal{E}}(\partial\Omega)$ yields a normal operator with spectrum contained in $\partial\Omega$. Therefore, $M_z$ on  the module $\mathcal{M}$ determined by $T_\varphi$  on $H^2_{\mathcal{E}}$ is a subnormal operator with its normal spectrum contained in $\partial\Omega$.  Hence, it is unitarily equivalent to a bundle shift (\cite{A-D}, Theorem 11). The fact that the multiplicity of the unitary representation  is finite follows from the fact that $\mathcal{M}$ has finite rank and hence so does the normal extension of $\mathcal{M}$.\qed

\vspace{0.3in}

If $H^2_\alpha(\Omega)$ is the module determined by a bundle shift and $\omega_0$ is a point in $\Omega$, then the submodule of $H^2_\alpha(\Omega)$ consisting of sections which vanish at $\omega_0$ is also a bundle shift. However, it need not be isometrically isomorphic to $H^2_\alpha(\Omega)$. But, the point $\omega_0$ can always be chosen so that it is. Thus the hypotheses of the theorem can be satisfied.

As was pointed out in the introduction, the Hardy spaces on finitely-connected domains are examples of \v Silov modules. Recall that a \v Silov module over a function algebra $A$ is a subnormal module for which the corresponding reductive or normal module containing it is actually over $C(\partial A)$, where $\partial A$ denotes the \v Silov boundary of $A$. It seems an interesting question as to whether these are the only finite rank, quasi-free Hilbert module containing a pure isometrically isomorphic submodule of finite codimension.

\vspace{0.6in}

\section{Essentially Reductive Case}

Now let us consider what we can say when the submodule $U\mathcal{R}$ has infinite codimension in $\mathcal{R}$. Let $\Omega$ be a domain in $\mathbb{C}^n$ and $\mathcal{R}$ be a quasi-free Hilbert module over $A(\Omega)$. Then the Hilbert space tensor product $\mathcal{R}\otimes H^2(\mathbb{D})$ is a quasi-free Hilbert module over $A(\Omega\times \mathbb{D})$ which clearly contains the pure isometrically isomorphic submodule $\mathcal{R}\otimes H^2_0(\mathbb{D})$. Hence, we can say little without some additional hypothesis  for $\Omega$ or $\mathcal{R}$ or both. One possibility would be to assume that $\Omega$ has no corners or is not a product. We will not pursue that direction here. Rather we make the additional assumption that $\mathcal{R}$ is essentially reductive. (It is an extremely interesting question as to whether essential reductivity is related to a lack of corners or not being a product.) The first result we obtain, seems at first, to be a little remarkable.

Recall that a Hilbert module $\clm$ is said to be essentially reductive \cite{D-P} if the operators $\{M_{\varphi}\}_{\varphi \in A(\Omega)}$ in 
$\cll(\clm)$ defined by
module multiplication  are all essentially normal, that is, the self-commutators 
$[M_{\varphi}^*, M_{\varphi}] = M_{\varphi}^* M_{\varphi} - M_{\varphi} 
M_{\varphi}^*$ are in the ideal of compact operators in $\clm$ for $\varphi$ in $A(\Omega)$.

\vspace{0.3in}
\begin{Theorem}\label{thm5}
Let $\mathcal{R}$ be an essentially reductive  Hilbert module over $A(\Omega)$ and $U$ be an isometric module map $U$  on $\mathcal{R}$ such that $\bigcap\limits^\infty_{k=0} U^k\mathcal{R} = (0)$. Then $\mathcal{R}$ is subnormal, that is, there exists a reductive Hilbert module $\mathcal{N}$ over $A(\Omega)$ with $\mathcal{R}$ as a submodule.
\end{Theorem}

\NI \textbf{Proof.} As before, there exists an isometric isomorphism $\Psi$ from $\mathcal{R}$ onto $H^2_{\mathcal{E}}(\mathbb{D})$ with $\mathcal{E} = \mathcal{R}\ominus U\mathcal{R}$ and $\varphi_1,\ldots, \varphi_n$ in $H^\infty_{\mathcal{L}(\mathcal{E})}(\mathbb{D})$ such that $\Psi$ is a $\mathbb{C}[\pmb{z}]$-module map relative to the module structure on $H^2_{\mathcal{E}}(\mathbb{D})$ defined so that $z_j\to T_{\varphi_j}$,  $j=1,2,\ldots, n$. We complete the proof by showing that the $n$-tuple $\{\varphi_1(e^{it}),\ldots, \varphi_n(e^{it})\}$ consists of commuting normal operators for $e^{it}$ a.e.\ on $\mathbb{T}$. Then $\mathcal{N}$ is $L^2_{\mathcal{E}}(\mathbb{T})$ with the module multiplication defined by $z_i\to L_{\varphi_i}$, where $L_{\varphi_i}$ denotes pointwise multiplication on $L^2_{\mathcal{E}}(\mathbb{T})$. Since the $\{\varphi_j(e^{it})\}^n_{j=1}$ are normal and commute, $L^2_{\mathcal{E}}(\mathbb{T})$ is a reductive Hilbert module.

The fact that $\mathcal{R}$ is essentially reductive implies that each $T_{\varphi_i}$ is essentially normal and hence that the cross-commutators $[T^*_{\varphi_i},T_{\varphi_j}]$ are  compact for $1\le i,j\le n$. We complete the proof by showing that $[T^*_{\varphi_i},T_{\varphi_j}]$ compact impiles that $[L^*_{\varphi_i}, L_{\varphi_j}] = 0$ on $L^2_{\mathcal{E}}(\mathbb{T})$.

Fix $f$ in $H^2_{\mathcal{E}}(\mathbb{D})$ and let $N$ be a positive integer. We observe that
\begin{equation}
\lim_{N\to\infty} \|(I-P) L^N_z L^*_{\varphi_i}L_{\varphi_j}f\| = 0\tag{*}
\end{equation}
and
\begin{equation}
\lim_{N\to\infty}\|(I-P) L^N_zL^*_{\varphi_i}f\| = 0,\tag{**}
\end{equation}
where $P$ is the projection of $L^2_{\mathcal{E}}(\mathbb{T})$ onto $H^2_{\mathcal{E}}(\mathbb{D})$. Therefore we have
\begin{align*}
\|[T^*_{\varphi_i}, T_{\varphi_j}] M^N_zf\| &= \|PL^*_{\varphi_i}PL_{\varphi_j}PL^N_zf - PL_{\varphi_j} PL^*_{\varphi_i}PL^N_zf\|\\
&= \|[L^N_zL^*{\varphi_i} L_{\varphi_j}f - (I-P) L^N_z L^*_{\varphi_i} L_{\varphi_j}f]\\
&\quad - [L_{\varphi_j} L^N_z L^*_{\varphi_i}f - L_{\varphi_j}(I-P) L^N_z L^*_{\varphi_i}f]\|.
\end{align*}
Using $(*)$ and $(**)$ we obtain
\begin{align*}
\lim_{N\to\infty} \|[T^*_{\varphi_i}, T_{\varphi_j}]L^N_zf\| &= \lim_{N\to\infty} \|(L^N_z L^*_{\varphi_i}L_{\varphi_j} - L_{\varphi_j}L^N_z L'_{\varphi_i})f\|\\
&= \lim_{N\to\infty} \|L^N_z[L^*_{\varphi_i}, L_{\varphi_j}]f\| = \|[L^*_{\varphi_i}, L_{\varphi_j}]f\|.
\end{align*}
Since $[T^*_{\varphi_i}, T_{\varphi_j}]$ is compact and the sequence $\{e^{iNt}f\}$ converges weakly to 0, we have $\lim\limits_{N\to\infty}\break \|[T^*_{\varphi_i},T_{\varphi_j}] e^{iNt}f\| = 0$. Therefore, $\|[L^*_{\varphi_i}, L_{\varphi_j}]f\| = 0$. Finally, the set of vectors $\{e^{-iNt}f\}\colon \ N\ge 0, f\in H^2_{\mathcal{E}}(\mathbb{D})\}$ is norm dense in $L^2_{\mathcal{E}}(\mathbb{T})$ and $\|[L^*_{\varphi_i},L_{\varphi_j}] e^{-iNt}f\| = \|[L^*_{\varphi_i}, L_{\varphi_j}]f\| = 0$. Therefore, $[L^*_{\varphi_i}, L_{\varphi_j}] = 0$ which completes the proof.\qed

\vspace{0.3in}

The following result is complementary to Theorem \ref{thm3}.

\vspace{0.3in}

\begin{Theorem} \label{reductive}
Let $\clm$ be an  essentially reductive, finite rank, quasi-free Hilbert module 
over $A(\mathbb{D})$. Let $U$ be a module isometry such that $\cap_{k = 0}^{\infty} 
U^k \clm = \{0\}$. Then $\clm$ is unitarily equivalent to 
$H^2_{\clf}(\mathbb{D})$ for some Hilbert space $\clf$ with dim$\clf$ = rank of 
$\clm$.
\end{Theorem}

\NI \textbf{Proof.} 
As before there is an isometrical isomorphism, $\Psi\colon \ H^2_{\mathcal{F}}(\mathbb{D})\to \mathcal{M}$ such that $U  = \Psi T_z\Psi^*$ and there exists $\varphi$ in $H^\infty_{\mathcal{L}(\mathcal{F})}(\mathbb{D})$ such that $M_z = \Psi T_\varphi \Psi^*$. Moreover, by Proposition \ref{pro3}, $\varphi(e^{it})$ is normal for $e^{it}$ in $\mathbb{T}$ a.e. Further, since $M_z$  is essentially normal and $M_z-\omega$ is Fredholm for $\omega$ in $\mathbb{D}$, it follows that $M_z$ is an essential unitary. Finally, this implies $T^*_\varphi T_\varphi-I = T_{\varphi^*\varphi-I}$ is compact and hence $\varphi^*(e^{it}) \varphi(e^{it}) = I$ a.e.\ or $\varphi$ is an inner function which completes the proof.\qed

\vspace{0.3in}

The only place in this proof in which the hypothesis that $\mathcal{M}$ is finite rank is needed is to conclude that $M_z-\omega$ is Fredholm for $\omega$ in $\mathbb{D}$. Thus it seems possible that the result would be true without that assumption.

If we consider the same question for $\mathbb{B}^n$ instead of $\mathbb{D}$, the result we obtain is that the Hilbert module is defined by a row isometry, which, of course, is not unique. In particular, one possibility is $H^2_{\mathcal{E}}(\mathbb{B}^n)$ but there are many others. For example, take $L^2_a(\mu)$ for any measure $\mu$ on $\partial\mathbb{B}^n$ for which $L^2_a(\mu) \ne L^2(\mu)$.

We give an application of these ideas for the $n$-shift space $H^2_n$ on the ball $\mathbb{B}^n$. In (\cite{GJX}, Corollary 5.5) Guo, Hu and Yu proved that two nested unitarily equivalent submodules of $H^2_n$ must be equal.
We provide a rather different proof when the larger module is $H^2_n$, which depends on analyzing the case when the module isometry is not pure. The  special properties of the $n$-shift space used are the fact that it is essentially reductive with essential module spectrum $\partial\mathbb{B}^n$ and on no submodule is it subnormal.

\begin{Corollary}
If $\mathcal{M}$ is a submodule of the $n$-shift space $H^2_n$ which is isometrically isomorphic to $H^2_n$, then $\mathcal{M} = H^2_n$. 
\end{Corollary}

\NI \textbf{Proof.} 
Again there exists an isometric module map $U$ on $H^2_n$ such that $\mathcal{M} = UH^2_n$. Let $U = U_s\oplus U_c \oplus U_p$ on $H^2_n = \mathcal{M}_s \oplus \mathcal{M}_c \oplus \mathcal{M}_p$ be the von Neumann--Wold decomposition (cf.\ \cite{S-N}) of $U$ on $H^2_n$ so that $U_s$ and $U_c$ are singular and absolutely continuous unitaries and $U_p$ is a pure isometry. Since there are no non-zero operators that intertwine $U_s$ with $U_c\oplus U_p$ in either direction, it follows that $\mathcal{M}_s$ is a reducing submodule of $H^2_n$ for the $\mathbb{C}[\pmb{z}]$ module action. However, an extension of the result in (\cite{C-D}, Corollary 1.11 or \cite{C-S}), shows that $H^2_n$ has no proper reducing submodules and  hence either $\mathcal{M}_s = (0)$ or $\mathcal{M}_s = H^2_n$. If $\mathcal{M}_s = H^2_n$, then $U$ is unitary and $\mathcal{M} = H^2_n$. Therefore, we may assume that $U$ on $H^2_n$ has no singular part.

Since the module $H^2_n$ is essentially reductive with essential module spectrum $\partial \mathbb{B}^n$, it follows that $I - \sum\limits^n_{i=1} M^*_{z_i}M_{z_i}$ is compact, where $M_{z_i}$ is the operator on $H^2_n$ defined to be module multiplication by $z_i$, $i=1,2,\ldots, n$. Since $\mathcal{M}_c$ is a submodule of $H^2_n$, it follows that $X = I - \sum\limits^n_{i=1}\widetilde M^*_{z_i} \widetilde M_{z_i}$ is compact, where $\widetilde M_{z_i}$ is the restriction of $M_{z_i}$ to $\mathcal{M}_c$ for $i=1,2,\ldots, n$. However, since $U$ is a module map, it follows that the absolutely continuous unitary operator commutes with $X$. Since $X$ is compact, it must be the zero operator and hence $\{\widetilde M_{z_i}\}$ is a row isometry. By the result of Athavale (\cite{A}, Proposition 2), the $n$-tuple is a jointly subnormal row isometry or $\mathcal{M}_c$ is a subnormal row isometry submodule of $H^2_n$. The following calculation uses the weights obtained by Arveson (\cite{A2}, Proposition 5.3) to show that $\mathcal{M}_c = (0)$.

Expand $f$ in $\mathcal{M}_c$, so that $f(\pmb{z}) = \sum\limits_{\pmb{\alpha}} a_{\pmb{\alpha}}z^{\pmb{\alpha}}$. Then
\begin{align*}
\sum^\infty_{k=0} \sum_{|\alpha|=k} |a_{\pmb{\alpha}}|^2 \|z^{\pmb{\alpha}}\|^2 &= \|f\|^2 =  \left\|\sum^n_{i=1} \widetilde M^\alpha_{z_i} \widetilde M_{z_i} f\right\| = \left\|\sum^n_{i=1} M^\alpha_{z_i}M_{z_i}f\right\| =\\
&= \sum^\infty_{k=0} \sum_{|\alpha|=k} \frac{1+k}{n+k} |a_{\pmb{\alpha}}|^2 \|z^{\pmb{\alpha}}\|^2,
\end{align*}
which implies $|a_{\pmb{\alpha}}| = 0$ for all multi-indices $\pmb{\alpha}$. As a consequence $\mathcal{M}_c = (0)$. Therefore, $U$ is pure and the result follows from the theorem.\qed

\vspace{0.3in}

If we consider the case of two nested submodules, then the approach used in this argument would require that they be essentially reductive which is not true for all submodules of $H^2_n$. However, there is another problem. One must eliminate somehow the possibility that $\mathcal{M}_s \ne (0)$. Unfortunately, one seems able only to conclude that $\mathcal{M}_s$ is a submodule of $H^2_n$ on which the module action yields a row contraction commuting with $U_s$ which is a row isometry modulo the compacts. But it is unclear how to show that this is not possible.

\vspace{0.5in}

\section{Submodules of Subnormal Modules}

\vspace{0.3in}

We conclude by considering when two submodules of a subnormal Hilbert module $\mathcal{M}$ over $A(\Omega)$ can be isometrically isomorphic. Of course, if $\mathcal{M}$ is \v Silov, then in many cases we have seen that there exists a proper isometrically isomorphic submodule. What about the non-\v Silov case?

If $\mu$ is the measure on $\text{clos } \mathbb{D}$ obtained from the sum of Lebesgue measure on $\partial\mathbb{D}$ and the unit mass at 0, then $L^2_a(\mu)$ is not a \v Silov module. However, it is easy to see that the cyclic submodules generated by $z$ and $z^2$, respectively, are isometrically isomorphic but distinct. A quick examination suggests the problem is that $\mu$ assigns positive measure to the intersection of a zero variety and $\mathbb{D}$. It turns out that if we exclude that possibility and $L^2(\nu)$ is not a \v Silov module, then distinct submodules can't be isometrically isomorphic. The proof takes several steps. 
\vspace{0.3in}
\begin{Lemma}
Let $\nu$ be a probability measure on $\text{clos } \Omega$ and $f$ and $g$ vectors in $L^2_a(\nu)$ so that the cyclic submodules of $L^2_a(\nu)$, $[f]$ and $[g]$, generated by $f$ and $g$, respectively, are isometrically isomorphic  with $f$ mapping to $g$. Then $|f| = |g|$ a.e. $\nu$. 
\end{Lemma}

\NI \textbf{Proof.} If the correspondence $Vf=g$ extends to an isometric module map, then $\langle \pmb{z}^{\pmb{\alpha}} f, \pmb{z}^{\pmb{\beta}}f\rangle_{L^2_a(\nu)}\break = \langle\pmb{z}^{\pmb{\alpha}}g, \pmb{z}^{\pmb{\beta}}g\rangle_{L^2_a(\nu)}$ for monomials $\pmb{z}^{\pmb{\alpha}}$ and $\pmb{z}^{\pmb{\beta}}$ in $\mathbb{C}[\pmb{z}]$. This implies that
\[
\int\limits_{\text{clos } \Omega} \pmb{z}^{\pmb{\alpha}} \bar{\pmb{z}}^{\pmb{\beta}} |f|^2 d\nu = \int\limits_{\text{clos } \Omega} \pmb{z}^{\pmb{\alpha}} \bar{\pmb{z}}^{\pmb{\beta}} |g|^2 d\nu\quad \text{for all monomials}\quad \pmb{z}^{\pmb{\alpha}}
\] 
and $\pmb{z}^{\pmb{\beta}}$. Since the linear span of the set $\{\pmb{z}^{\pmb{\alpha}} \bar{\pmb{z}}^{\pmb{\beta}}\}$ forms a self-adjoint algebra which separates the points of clos $\Omega$, it follows that the two measures $|f|^2~d\nu$ and $|g|^2~d\nu$ are equal or that $|f| = |g|$ a.e. $\nu$.\qed

\vspace{0.3in}

Although it might be possible to avoid it, we consider measures for which point evaluation on $\Omega$ is bounded.

\vspace{0.3in}
\begin{Theorem}\label{thm1}
Let $\nu$ be a probability measure on $\text{\rm clos } \Omega$ such that point evaluation is bounded on $L^2_a(\Omega)$ with closed support properly containing $\partial\Omega$ but such that $\nu(X) = 0$ for $X$ the intersection of $\text{\rm clos } \Omega$ with a zero variety. If $\mathcal{M}_1$ and $\mathcal{M}_2$ are isometrically isomorphic submodules of $L^2_a(\nu)$, then $\mathcal{M}_1 = \mathcal{M}_2$.
\end{Theorem}

\NI\textbf{Proof.} Let $V$ be an isometric module map from $\mathcal{M}_1$ onto $\mathcal{M}_2$. For $0\ne f$ in $\mathcal{M}_1$, let $g = Vf$. By the previous lemma, we have $|f| = |g|$ a.e.\ $\nu$. Since $\partial\Omega$ is contained in the closed support of $\nu$, it follows that $|f(\pmb{\omega})| = |g(\pmb{\omega})|$ for $\pmb{\omega}$ in $\partial\Omega$. Since point evaluation is bounded on $L^2_a(\Omega), f$ and $g$ are holomorphic on $\Omega$. If $X = \{\pmb{\omega} \in\Omega \mid f(\pmb{\omega}) = 0\}$,  then $\nu(X) = 0$ which implies that $\nu(\Omega\backslash X)>0$. If we set $h(\pmb{\omega}) = \frac{g(\pmb{\omega})}{f(\pmb{\omega})}$ for $\pmb{\omega}$ in $\Omega\backslash X$, then $\sup\limits_{\pmb{\omega}\in\Omega\backslash X} |h(\pmb{\omega})|\le 1$. Since there is $\pmb{\omega}_0$ in the support of $\nu$ in $\Omega\backslash X$ such that $|h(\pmb{\omega}_0)| = 1$, we have $|h(\pmb{\omega})| \equiv 1$ on $\Omega\backslash X$. (Here we are using the fact that $\Omega \backslash X$ is connected.) Thus there is a constant $e^{i\theta}$ such that $f=e^{i\theta}g$ on $\Omega$.

Since this holds for every $f$ in $\mathcal{M}_1$, by considering $f_1,f_2$ and $f_1+f_2$, we see that $Vf=e^{i\theta}f$ for all $f$ in $\mathcal{M}_1$ and hence $\mathcal{M}_1 = \mathcal{M}_2$. \qed

\vspace{0.3in}

This result contains the results of Richter \cite{Rich}, Putinar \cite{P}, and Guo--Hu--Xu \cite{GJX} mentioned earlier since area measure on $\mathbb{D}$ or volume measure on $\Omega$ satisfies the hypotheses of the theorem. However, so do the measures for the weighted Bergman spaces on $\mathbb{D}$ or weighted volume measure on any domain $\Omega$.

\begin{Corollary}
If $\Omega$ is a bounded domain in $\mathbb{C}^n$ and $\mathcal{M}_1$ and $\mathcal{M}_2$ are isometrically isomorphic submodules of $L^2_a(\Omega)$, then $\mathcal{M}_1 = \mathcal{M}_2$.
\end{Corollary}

If $\mathcal{R}$ is a subnormal Hilbert module over $A(\Omega)$ of finite multiplicity greater than one, then the conclusion of the previous results doesn't follow. One possible substitute result would be the existence of a unitary module map on the normal module which extends it which takes one module to the other. That is not quite right but perhaps something like that is. We leave the formulation of such a result as an open problem.

\end{document}